\documentclass{amsart}

\usepackage{graphics}
\usepackage{epsfig}
\usepackage{amssymb,amsmath,mathrsfs}
\usepackage{amsthm}
\usepackage{amsfonts}
\usepackage[english]{babel}
\usepackage{blindtext}
\usepackage{color}
\usepackage{hyperref}
\usepackage{algorithmic}
\usepackage{algorithm}
\usepackage[table]{xcolor}
\usepackage{ragged2e}
\usepackage{dcolumn}

\newtheorem{theorem}{Theorem}[section]

\theoremstyle{definition}

\newtheorem{example}[theorem]{Example}

\theoremstyle{remark}

\newcommand{\dom}{\mathop{\rm dom}\nolimits}

\numberwithin{equation}{section}

\begin{document}

\title[Newton-GSOR method for  Unconstrained Optimization Problems]{Newton-GSOR method for solving Large-Scale Unconstrained Optimization Problems}


\author{Santoshi Subhalaxmi Ray}
\address{Department of Mathematics, National Institute of Technology Meghalaya, shillong, meghalaya, 793003}
\curraddr{}
\email{s21ma008@nitm.ac.in}
\thanks{}

\author{Manideepa Saha}
\address{Department of Mathematics, National Institute of Technology Meghalaya, shillong, meghalaya, 793003}
\curraddr{}
\email{manideepa.saha@nitm.ac.in}
\thanks{}


\date{}

\dedicatory{}


\begin{abstract} Unconstrained convex optimization problems have enormous applications in various field of science and engineering. Different iterative methods are available in literature to solve such problem, and Newton method is among the oldest and simplest one. Due to slow convergence rate of Newton's methods, many research have been carried out to modify the Newton's method for faster convergence rate.  In 2019, Ghazali et al. modified Newton’s method and proposed Netwon-SOR method, which is a combination of Newton method with SOR iterative method to solve a linear system. In this paper, we propose a modification of Newton-SOR method by modifying SOR method to generalized SOR method. Numerical experiments are carried out to check the efficiently of the proposed method.

\end{abstract}

\maketitle

\section{Introduction}\label{sec1}
 
Consider a large scale unconstrained optimization problem which is of the form 
\begin{equation}\label{eq1}
\min f(x)
\end{equation}
where $f:\mathbb{R}^n\to \mathbb{R}$,is the objective function, which is convex and twice continuously differentiable. Unconstrained convex optimization problem doesn't include neither any inequality constraints nor any equality constraints. Unconstrained convex optimization problems have enormous applications in  various field of science and engineering. These problems appears in  image processing, electric power systems,  DNA reproduction system etc. Furthermore, constrained optimization problems may also solved as a sequence of unconstrained optimization problem. We refer to~\cite{NLL02,LL95,MM18,ZL17}  for further applications.

Various numerical methods are developed to solve unconstrained optimization problem~\eqref{eq1}, and these methods are categories in two types:  gradient methods and direct search (non gradient methods) methods. The direct search methods doesn't involve the partial derivative of the objective functions.  Rosenbrock’s method\rm\cite{R60},  Hooke and Jeeve’s method \rm\cite{YB94}, Nelder-Mead simplex method \rm\cite{DW16} are some popular the direct search methods. In contrast to the direct search methods, gradient methods involves the partial derivative. Steepest descent method\rm\cite{NSMHS18}, conjugate gradient method\rm\cite{ML16}, Newton’s method ,Quasi-Newton method\rm\cite{AAD15} etc., are some  existing  gradient method. Newton's method is one of the oldest best-known techniques because of it's quadratic convergence.  The was proposed by Newton in 1669. This method is also considered as the fundamental for the most effective methods both constrained and unconstrained optimization problems. For example, in convex optimization, polynomial time interior point methods are developed from the Newton’s method.  We refer to the survey paper~\cite{P07} for further details.

As $f$ is a convex differentiable function, a necessary and sufficient condition for the optimum solution $\tilde{x}$ of~\eqref{eq1} is
\begin{equation}\label{eqn2}
\nabla f(\tilde{x})=0
\end{equation}
where $\nabla f(x)$ represents the gradient vector of $f$. Thus solving~\eqref{eq1} is equivalent to find the solution of~\eqref{eqn2}. It may be noted that ~\eqref{eqn2} is collection of $n$ equations in $n$ unknown, usually such problems can be solved by iterative method. A line search method to solve~\eqref{eqn2} and hence of~\eqref{eq1} consisting of the iterative  step
\begin{equation}\label{eqn3}
x^{(k+1)}=x^{(k)}+\alpha_{k}d_k
\end{equation}
where $\alpha_k$ is positive integer that represents the step length and $d_k$ is the search direction. Various choices of $\alpha_k$ and $d_k$ defines various line search methods.  Descent method to solve~\eqref{eq1} is among one such line search method in which the search direction takes the form
\begin{equation}\label{eqn4}
d_k=-C_k(\nabla (x^{(k)})
\end{equation}
for some nonsingular symmetric matrix $C_k$. If $C_k=-I$, then  equation~\eqref{eqn3} is known as the steepest descent method. By taking $-C_k$ as the exact Hessian matrix and approximate Hessian matrix, the method~\eqref{eqn2} is  respectively, known as Newton and quasi-newton method.  The search direction in Newton's method assured that
\begin{equation}\label{eqn5}
f (x^{(k+1)})<f (x^{(k)})
\end{equation}
If the initial approximation $x^{(0)}$ chosen in Newton's  method is closer to the optimal solution $\tilde{x}$, then Newton iteration converges rapidly to $\tilde{x}$. However, the converges rate may be slower if $x^{(0)}$ is away from $\tilde{x}$. 
Due to the numerous applications and simplicity of Newton's method, researchers  developing methods based on Newton's method to make the convergence faster and to ensure the convergence of the method.  In this paper, we propose a combination of Newton's method with generalised SOR method for faster  convergence of Newton's method. In~\cite{GLL86}, authors proposed  a nonmonotone  step size selection rule in Newton's method and numerical illustration showed that faster convergence of the method. In~\cite{S20}, author  presented a updation of Newton's method by combining Newton's direction and the steepest descent direction in each iterative step. The method assured the convergence for any initial approximation$x^{(0)}$. Recently, in 2019~\cite{GSDG19},  Ghazali proposed a new algorithm to solve~\eqref{eqn2} which is a combination of Newton's method and successive over relaxation( SOR ) method, (called as Newton-SOR method) and analyzed the performance of the method. Motivated by their work, in this paper we propose a modification of Newton-SOR method by replacing SOR method with generalized method and discuss the performance of the proposed algorithm. We refer the proposed algorithm as Newton-GSOR method.

The paper is organized as follows: In Section~2, we discuss the formulation of the method and propose the algorithm. In Section~3, we consider various numerical illustrations to analyze the effectiveness of the proposed algorithm. At last, we ended with a overall conclusion in Section~4.

\section{Formulation of Newton-GSOR method} In this section we discuss formulation of  the steps of Newton-GSOR method.  We begin with an well-known result of calculus, which is given below:

\begin{theorem}\label{thm1}\em\cite{BV04} A necessary and sufficient criterion for the optimality of ~\eqref{eq1} at a point $x^*$ is 
\[\nabla f(x^*)=0\]
where $\nabla f(x^*)$ represents the gradient vector of $f$ at $x^*$.
\end{theorem}
We revisit the steps for Newton's method. To minimize the objective function $f(x)$ at the point $x^{(k)}$, we  first show that steps to formulate the Newton iteration, which is obtained  by neglecting higher order terms of Taylor series expansion:
\begin{equation}\label{eq2}
f(x)\approx f(x^{(k)})+[\nabla f(x^{(k)})]^T(x-x^{(k)})+\frac{1}{2}(x-x^{(k)})^T H(f(x^{(k)}))(x-x^{(k)})
\end{equation} 
where $H(f(x))$ represents the Hessian matrix of $f$ at the point $x$. We simply write $H(x)$ for  $H(f(x))$. Assuming that $x^{(k)}$ is an optimal solution~of~\eqref{eq1} and using Theorem~\ref{thm1} we get
\begin{equation}
x=x^{(k)}-\left[H\left(x^{(k)}\right)\right]^{-1}\nabla f(x^{(k)})~,
\end{equation}
This suggests us to choose the next step of iteration as follows:
\begin{equation}\label{eq3}
x^{(k+1)}=x^{(k)}-\left[H\left(x^{(k)}\right)\right]^{-1}\nabla f(x^{(k)})
\end{equation}
Since the descent direction of the  Newton method is $d^{(k)}=x^{(k+1)}-x^{(k)}$ which can be obtained as follows:
\begin{equation*}
d^{(k)}= -\left[ H\left(x^{(k)}\right)\right]^{-1}\nabla f\left(x^{(k)}\right)\implies H\left(x^{(k)}\right)  d^{(k)}=-\nabla f\left(x^{(k)}\right)
\end{equation*}
This shows that the descent direction can be obtained  by solving the following linear system
\begin{equation}\label{eq4}
Hd=\hat{f}~,
\end{equation}
where $H$ is the Hessian matrix $H(x^{(k)})$, $d=[d_1~ d_2~ d_3~\dots~d_n]^T$, and $\hat{f}=-\nabla f(x^{(k)})$ at the $k$-th step.  Thus the Newton's algorithm is  given by:

\begin{algorithm}[ht]
\caption{\label{alg1}}
	\begin{flushleft} 
	\begin{itemize}
\item[1.]Choose a starting point $x^{(0)}\in \dom(f)$.
\item[2.] Until stopping criterion is satisfied, do
\begin{itemize}
 \item [i).] Calculate  $d_k:=-H(x^{(k)})^{-1}\nabla f(x^{(k)})=-H^{-1}\hat{f}$. If $d_k=0$, then STOP
 \item[ii).] Choose a step size $\alpha_k=1$. 
\item[iii).] Update. $x^{(k+1)}:= x^{(k)} +\alpha_{k}x^{(k)}$.  	
\end{itemize}
	\end{itemize}
	\end{flushleft} 
\end{algorithm}
For large value of $n$, calculating $H(x^{(k)})^{-1}$ and hence of $d_{k}$ takes a large memory of the  computer. So, we find~$d_k$ in Step 2(i) of the above algorithm by generalized SOR method proposed in~\cite{SC20}. We  now discuss generalized SOR (GSOR) method introduced in~\cite{SC20}. In the GSOR method the coefficient matrix  $H_k$ is decomposed as
\begin{equation}
H_k=T_m-E_m-F_m,
\end{equation}
where $T_m=(t_{ij})$ is a banded matrix with bandwidth $2m+1$ defined as
\begin{equation}\label{band-mat}
t_{ij} = \left\{
        \begin{array}{ll}
            a_{ij}, & |i-j| \leq m \\
              0 ,& \text{ otherwise }
        \end{array}
    \right.
    \end{equation}
and $E_m$ and $F_m$ are the strict lower part and upper part of $H$.  More specifically, the  matrices  $T_m,~E_m$ and $F_m$ are defined as following

\begin{eqnarray}
T_{m}=\left[\begin{array}{rrrr}
 a_{1,1} & \dots & a_{1,m+1} & 0 \\
\vdots  & \ddots & \ddots & \\
a_{m+1,1 } &  & & a_{n-m,n}\\
 & \ddots&\ddots & \\ 
0 & & a_{n,n-m} & a_{n,n}
\end{array}\right],\\\nonumber
E_{m}=\left[\begin{array}{rrr}
  0 & \dots & 0 \\
 -a_{m+2,1}   \\
\vdots &\ddots&\vdots  \\
-a_{n,1} & \dots & -a_{n-m-1,n}
\end{array}\right] \nonumber
\end{eqnarray}

\begin{equation*}\label{eqn2.2}
F_{m}=\left[\begin{array}{rrrl}
0 & -a_{1,m+2}  &\ldots -a_{1,n} \\
\vdots &\ddots   &\vdots \\
0 & \ldots & -a_{n-m-1,n}\\ & &
\end{array}\right]
\end{equation*}
The iterative formula of Generalized Jacobi (GJ)method can be stated as 
\begin{equation}\label{eq5}
d^{(k+1)}=T_m^{-1}(E_m+F_m)~d^{(k)}+T_m^{-1}\hat{f}
\end{equation}
Here we are considering the splitting as $H=M-N$. For Generalized Jacobi method, $M=T_m$ and $N=E_m+F_m$.\\

Similarly, the iterative formula of Generalized Gauss-Seidal (GGS) method is given by
\begin{equation}\label{eq6}
d^{(k+1)}=(T_m-E_m)^{-1}F_m~d^{(k)}+(T_m-E_m)^{-1}\hat{f},
\end{equation}
considering the splitting as $M=T_m-E_m$ and  $N=F_m$.\\

If we take $m=0$, then \eqref{eq5} and \eqref{eq6} reduces to Jacobi and Gauss-Seidel, respectively. The iterative scheme for the GSOR method to solve~\eqref{eq4} is given by
\begin{equation}\label{eq7}
d^{(k+1)}=(T_m-\omega E_m)^{-1}(\omega F_m+(1-\omega) T_m)~d^{(k)}+\omega (T_m-\omega E_m)^{-1}\hat{f}
\end{equation}
It may be noticed that in case $\omega=1$, GSOR method will reduced to GGS method. We now propose the Newton-GSOR algorithm for~\eqref{eq1} 

\newpage
\begin{algorithm}[ht]
\caption{\label{alg2}}
	\begin{flushleft} 
	\begin{itemize}

\item[1.] Initialize $f(x) \in \mathbb{R},~ x^{(0)}\in \mathbb{R}^n,~\omega\in(1,2]$ and $\epsilon_1,~\epsilon_2$ 
\item[2.] Set the optimal value of $\omega$. 
\item[3.]  Until stopping criterion $\lVert f(x^{(k)}) \rVert < \epsilon_1~$ is satisfied, do
\begin{itemize}
\item[a.] Assign $d^{(0)}=0$.
\item[b.] Compute $\hat{f}:=\Delta f(x^{(k)}),~H=H(f(x^{(k)}))$
\item[c.] Decompose $H=T_m-E_m,F_m$
\item[d.] Compute the approximate value of $d^{(k+1)}$ as follows:
\begin{center}
$d^{(k+1)}=(T_m-\omega E_m)^{-1}(\omega F_m+(1-\omega) T_m)~d^{(k)}+\omega (T_m-\omega E_m)^{-1}\hat{f}$
\end{center}
 \item[e.] Check the convergence, $\lVert d^{(k+1)}-d^{(k)} \rVert < \epsilon_2~$. If it satisfied,  go to next step~(e), otherwise go back to step(b).
 \item[e.] Calculate $x^{(k+1)}=x^{(k)}+d^{(k)}$ .
\end{itemize}

\end{itemize}
	\end{flushleft} 
\end{algorithm}

\section{Numerical Experiments} In this section numerical experiments are carried out to illustrate the efficiency of Newton-GSOR method with Newton-SOR method.  The computations are carried out in MATLAB(2022b)  on a machine with Processor 3.30 GHz AMD Ryzen 5 5600H and CPU 8 GB . The computations are rounding to four digits  with the threshold $\epsilon_1=10^{-6}$ and $\epsilon_2=10^{-8}$. 

\begin{example}\rm\cite{A08} The LIARWHD Function is defined as
\begin{equation}
f(x)=\sum_{i=1}^{n}4(x_i^2-x_1)^2+\sum_{i=1}^{n}(x_i-1)^2
\end{equation}
This function achieves global minimum $f^*=0$ at $x^*=(1,1,\ldots,1)$. Applying  Newton-GSOR method for various values of $n$ and $m$, a comparison table  in terms of iteration count (IC) and time (T) with Newton-SOR method is given below. We also  provide  the comparison with Newton-GJ and Newton-GGS. For this numerical experiment, the initial guess $x^{(0)}$ is considered as $x^{(0)}=(4,4,\ldots,4)$.

\begin{table}[ht]
\begin{center}
\begin{tabular}{||p{1cm}|p{1cm}|p{1cm}|p{1.5cm}|p{1.4cm}|p{1cm}|p{1cm}|c||}
	\hline
	\multicolumn{4}{|c|}{Newton-SOR}&\multicolumn{4}{c|}{Newton-GSOR} \\ \hline
	n & Outer IC &  Inner IC & { T (in sec)} & m & Outer IC &  Inner IC& T (in sec)\\ \hline
	20 & 11 & 218 & 97.17 & 15 & 11 & 96 & 92.517\\ \hline
	30 & 11 & 244 & 303.579 & 25 & 11 & 90 & 293.079\\ \hline
	50 & 11 & 275 & 1298.399 & 45 & 11 & 88 & 1291.239\\ \hline
	100 & -- & -- & -- & 95 & 11 & 93 & 9786.021\\ \hline
\end{tabular}
\end{center}
\vspace{2mm}
\caption{Comparison table  between NSOR and NGSOR }
\end{table}

\newpage
\begin{table}[ht]
\begin{center}
\begin{tabular}{||p{1cm}|p{1cm}|p{1cm}|p{1cm}|p{1.5cm}|p{1.5cm}|p{1.5cm}|p{1.2cm}||}
	\hline
	 \multicolumn{5}{|c|}{Newton-GGS} &\multicolumn{3}{|c|}{Newton-GJ} \\ \hline
	n & m &  Outer IC &  Inner IC & T (in sec) &  Outer IC &  Inner IC & T (in sec) \\ \hline
	20  & 15  & 11 & 115 & 92.792 & 11 & 204 & 93.831 \\ \hline
	30  & 25  & 11 & 105 & 295.85 & 11 & 183 & 299.692 \\ \hline
	50  & 45  & 11 & 92 & 1292.541 & 11 & 161 & 1300.143 \\ \hline
   100  & 95  & 11 & 80 & 9808.256 & -- & -- & -- \\ \hline
\end{tabular}
\end{center}
\vspace{2mm}
\caption{Comparison table  between  NGGS and NGJ }
\end{table}
Changing the initial guess to $x^{(0)}=(1.5,1.5,\ldots,1.5)$, we obtained the following data.  However the experiment is done only for $n=20$ and $n=30$.

\begin{table}[ht]
\begin{center}
\begin{tabular}{||p{1cm}|p{1cm}|p{1cm}|p{1.5cm}|p{1.4cm}|p{1cm}|p{1cm}|c||}
	\hline
	\multicolumn{4}{|c|}{Newton-SOR}&\multicolumn{4}{c|}{Newton-GSOR} \\ \hline
	n & Outer IC &  Inner IC & { T (in sec)} & m & Outer IC &  Inner IC& T (in sec)\\ \hline
	20 & 8 & 143 & 66.528 & 15 & 8 & 64 & 65.491\\ 
		 &  &  &  & 17 & 8 & 56 & 60.446\\ \hline
	30 & 8 & 164 & 210.146 & 25 & 8 & 57 & 203.346\\ 
	 &  &  &  & 27 & 8 & 56 & 203.243\\ \hline
	
\end{tabular}
\end{center}
\vspace{2mm}
\caption{Comparison table  between NSOR and NGSOR }
\end{table}

\begin{table}[ht]
\begin{center}
\begin{tabular}{||p{1cm}|p{1cm}|p{1cm}|p{1cm}|p{1.5cm}|p{1.5cm}|p{1.5cm}|p{1.2cm}||}
	\hline
	 \multicolumn{5}{|c|}{Newton-GGS} &\multicolumn{3}{|c|}{Newton-GJ} \\ \hline
	n & m &  Outer IC &  Inner IC & T (in sec) &  Outer IC &  Inner IC & T (in sec) \\ \hline
	20  & 15  & 8 & 93 & 65.51 & 8 & 168 & 66.684 \\ 
	  & 17  & 8 & 68 & 69.695 & 8 & 119 & 67.486 \\ \hline
	30  & 25  & 8 & 83 & 203.489 & 8 & 148 & 205.91 \\
	  & 27  & 8 & 60 & 209.07 & 8 & 107 & 205.166 \\ \hline
	
\end{tabular}
\end{center}
\vspace{2mm}
\caption{Comparison table  between  NGGS and NGJ }
\end{table}

\end{example}

\begin{example}\rm\cite{A08} The DIAG-AUP1 function is defined as
\begin{equation}
f(x)=\sum_{i=1}^{n}4(x_i^2-x_1)^2+\sum_{i=1}^{n}(x_i^2-1)^2
\end{equation}
 This function also achieves global minimum $f^*=0$ at $x^*=(1,1,\ldots,1)$. Applying  Newton-GSOR method for various values of $n$ and $m$ and with the initial guess $x^{(0)}=(4,4,\ldots,4)$, a comparison table with Newton-SOR method is given below:
 \begin{table}[ht]
\begin{center}
\begin{tabular}{||p{1cm}|p{1cm}|p{1cm}|p{1.5cm}|p{1.4cm}|p{1cm}|p{1cm}|c||}
	\hline
	\multicolumn{4}{|p{5.5cm}|}{Newton-SOR}&\multicolumn{4}{p{5.5cm}|}{Newton-GSOR} \\ \hline
	n & Outer IC &  Inner IC & { T (in sec)} & m & Outer IC &  Inner IC& T (in sec)\\ \hline
	20 & 12 & 404 & 258.652 & 15 & 12 & 165 & 121.834\\ \hline
	30 & 12 & 493 & 814.757 & 25 & 12 & 150 & 813.269\\ \hline
	50 & 12 & 734 & 2953.88 & 45 & 12 & 131 & 1407.675\\ \hline
	100 & -- & -- & -- & 95 & 13 & 111 & 30952.604\\ \hline
\end{tabular}
\end{center}
\vspace{2mm}
\caption{Comparison table  between NSOR and NGSOR }
\end{table}

\begin{table}[ht]
\begin{center}
\begin{tabular}{||p{1cm}|p{1cm}|p{1cm}|p{1cm}|p{1.5cm}|p{1.5cm}|p{1.5cm}|p{1.4cm}||}
	\hline
	 \multicolumn{5}{|p{4.5cm}|}{Newton-GGS} &\multicolumn{3}{|p{3.5cm}|}{Newton-GJ} \\ \hline
	n & m &  Outer IC &  Inner IC & T (in sec) &  Outer IC &  Inner IC & T (in sec) \\ \hline
	20  & 15  & 12 & 200 & 124.553 & 12 & 359 & 271.234 \\ \hline
	30  & 25  & 12 & 181 & 806.066 & 12 & 328 & 757.854 \\ \hline
	50  & 45  & 12 & 157 & 1211.499 & 12 & 284 & 1404.33 \\ \hline
   100  & 95  & 13 & 133 & 19781.485 & 13 & 240 & 24269.481 \\ \hline
\end{tabular}
\end{center}
\vspace{2mm}
\caption{Comparison table  between  NGGS and NGJ }
\end{table}

Following data table is obtained for the initial guess $x^{(0)}=(1.5,1.5,\ldots,1.5)$, and for different values of the bandwidth.

\begin{table}[ht]
\begin{center}
\begin{tabular}{||p{1cm}|p{1cm}|p{1cm}|p{1.5cm}|p{1.4cm}|p{1cm}|p{1cm}|c||}
	\hline
	\multicolumn{4}{|c|}{Newton-SOR}&\multicolumn{4}{c|}{Newton-GSOR} \\ \hline
	n & Outer IC &  Inner IC & { T (in sec)} & m & Outer IC &  Inner IC& T (in sec)\\ \hline
	20 & 8 & 360 & 85.089 & 15 & 8 & 127 & 66.921\\ 
	 &  &  &  & 17 & 8 & 80 & 66.205\\ \hline
	30 & 9 & 448 & 249.123 & 25 & 8 & 113 & 206.105\\ 
	 &  &  &  & 27 & 8 & 74 & 204.874\\ \hline
	
\end{tabular}
\end{center}
\vspace{2mm}
\caption{Comparison table  between NSOR and NGSOR }
\end{table}

\begin{table}[ht]
\begin{center}
\begin{tabular}{||p{1cm}|p{1cm}|p{1cm}|p{1cm}|p{1.5cm}|p{1.5cm}|p{1.5cm}|p{1.2cm}||}
	\hline
	 \multicolumn{5}{|c|}{Newton-GGS} &\multicolumn{3}{|c|}{Newton-GJ} \\ \hline
	n & m &  Outer IC &  Inner IC & T (in sec) &  Outer IC &  Inner IC & T (in sec) \\ \hline
	20  & 15  & 8 & 169 & 68.331 & 8 & 315 & 69.674 \\ 
	  & 17  & 8 & 112 & 66.687 & 8 & 203 & 66.961 \\ \hline
	30  & 25  & 8 & 150 & 206.827 & 8 & 283 & 214.447 \\ 
	  & 27  & 8 & 102 & 202.704 & 8 & 187 & 208.084 \\ \hline
	
\end{tabular}
\end{center}
\vspace{2mm}
\caption{Comparison table  between  NGGS and NGJ }
\end{table}
\end{example}

\newpage
From the above two examples we have observed that Newton-GSOR methods takes less number of iterations in compare to the Newton-SOR method. Further, we notice that no. of iteration in Newton-GSOR method is less with increase bandwith.

\section{Conclusion} In this paper we proposed a modification of Newton-SOR method, and named it as Newton-GSOR method. The proposed method is a combination of Newton's method and generalized SOR method~\cite{SC20}. We provided an algorithm for the proposed method. We considered two test functions from \cite{A08} for numerical experiments. From the numerical illustration, we observed that Newton-GSOR method is faster  ( both in iteration count and in time) than Newton-SOR, Newton-GJ and Newton GGS method. Furthermore, it may be noticed that Newton-GSOR method takes almost half iteration ( time) than that of Newton-SOR method. Further, we have noticed that if bandwidth $m$ increases, then the convergence rate also  increases.

\end{document}